\documentclass[reqno,a4paper,12pt]{amsart} 

\usepackage{amsmath,amscd,amsfonts,amssymb}
\usepackage{mathrsfs,dsfont}
\usepackage{color}
\usepackage{mathtools}
\usepackage[hidelinks]{hyperref}
\hypersetup{colorlinks=true,linktocpage=true,
	linkcolor=blue,
	filecolor=blue,
	citecolor = blue,
	urlcolor=red}
\usepackage{tikz-cd}

\usepackage{tikz}
\usetikzlibrary{decorations.pathreplacing}

\numberwithin{equation}{section}
\numberwithin{figure}{section}

\def\a{\alpha}

\def\b{\beta}
\def\e{\varepsilon}
\def\g{\gamma}

\def\d{\Delta}
\def\n{\nabla w}
\def\n2{\nabla^2 w}

\def\R{\mathbb{R}}

\def\M{\mathcal{M}}

\renewcommand\le{\leqslant}
\renewcommand\ge{\geqslant}
\renewcommand\leq{\leqslant}
\renewcommand\geq{\geqslant}

\theoremstyle{plain}
\newtheorem{thm}{Theorem}[section]
\newtheorem*{thm*}{Theorem}
\newtheorem{theorem}[thm]{Theorem}
\newtheorem{lem}[thm]{Lemma}
\newtheorem{lemma}[thm]{Lemma}
\newtheorem{corollary}[thm]{Corollary}

\newtheorem{problem}[thm]{Problem}

\newtheorem*{claim*}{Claim}

\theoremstyle{definition}

\newtheorem*{definition*}{Definition}
\newtheorem*{remarks*}{Remarks}
\newtheorem*{remark*}{Remark}
\newtheorem{remark}[thm]{Remark}
\newtheorem{example}[thm]{Example}

\newenvironment{enumerate-math}
{\begin{enumerate}
		\addtolength{\itemsep}{5pt}
		}
	{\end{enumerate}}

\newenvironment{enumerate-text}
{\begin{enumerate}
		\addtolength{\itemsep}{5pt}
		}
	{\end{enumerate}}

\makeatother
\makeatletter
\@namedef{subjclassname@2020}{\textup{2020} Mathematics Subject Classification}
\makeatother

\begin{document}

\title[gradient estimates and universal bounds]
{logarithmic gradient estimates and\\ universal bounds for semilinear\\ elliptic equations revisited}
	
	\author{Zhihao Lu}
	\address{School of Mathematics and Statistics, Jiangsu Normal University, Xuzhou 221116, P. R. China}
	\email{zhihaolu@jsnu.edu.cn}

	\date{}
	
	
	\begin{abstract} 
		We derive the complete and optimal Cheng--Yau  gradient estimates and universal bounds for subcritical semilinear elliptic equations on complete Riemannian manifolds with (Bakry-\'{E}mery) Ricci curvature bounded below. This resolves a long-standing fundamental problem raised in Li’s work \cite{LJY}. As a corollary, this provides a new proof of the classical Gidas-Spruck  Liouville theorem. The Harnack inequality is also obtained.
	\end{abstract}
	
	\maketitle
	
	\section{Introduction}

	\subsection{Motivation and results}
	The study of harmonic functions on Riemannian manifolds has been one of the central subjects in geometric analysis. In their seminal work \cite{CY}, Cheng and Yau derived the following gradient estimate for positive harmonic functions on Riemannian manifolds:

	\begin{thm*}[Cheng--Yau]
		Let \({\mathcal{M}}^{n}\) be an \(n\)-dimensional complete Riemannian manifold with Ric \(\geq- {Kg}\) , where \(K \geq  0\) is a constant. Suppose that \(u\) is a positive harmonic function on a geodesic ball \(B\left( {{x}_{0},R}\right)\). Then
		\begin{equation}\label{111}
			\mathop{\sup }\limits_{{B\left( {{x}_{0},R/2}\right) }}\frac{{\left| \nabla u\right| }^{2}}{{u}^{2}} \leq  C\left( n\right) \left( {K + \frac{1}{{R}^{2}}}\right), 
		\end{equation}
		where \(C\left( n\right)\) is a constant depending only on \(n\).
	\end{thm*}
	
An important feature of Cheng-Yau's estimate is that the right-hand side of \eqref{111} depends only on \(n,K\) and \(R\) , and it does not depend on the lower bound of the injectivity radius or a global coordinate system. Notice that the important Harnack inequality and Liouville theorem are also obtained from their gradient estimate and the dependence of the Harnack constant on the lower bound of the curvature is clear. Afterwards, Cheng-Yau's approach deeply influenced many important results in linear geometric analysis such as the eigenvalue estimates of a manifold in Li-Yau \cite{LY0} and Zhong-Yang \cite{ZY}, Li-Yau's differential Harnack inequality for heat equation in \cite{LY}, etc. We refer to the classical books of Schoen-Yau \cite{SZ} and Peter Li \cite{LPE} for an overview of the subject.

Since 1980 and thereafter, many geometric analysts have considered the properties of nonlinear equations on manifolds following the same philosophy. The following Lane-Emden equation is a fundamental model of semilinear elliptic equation.
\begin{equation}\label{LE}
	{\Delta u} + {u}^{\alpha } = 0,
\end{equation}
where \(\alpha  \in  \mathbb{R}\). In 1981, Gidas and Spruck obtained the following celebrated Liouville theorem for subcritical equation \eqref{LE} in the work \cite{GS}.
	
\begin{thm*}[Gidas--Spruck]
	Let \(\left( {{\mathcal{M}}^{n},g}\right) \left( {n \geq  3}\right)\) be an \(n\)-dimensional complete Riemannian manifold with non-negative Ricci curvature. Let $u$ be a non-negative \({C}^{2}\) solution of \eqref{LE} on \({\mathcal{M}}^{n}\) with \(\alpha  \in  \left( {1,\frac{n + 2}{n - 2}}\right)\) , then \(u \equiv  0\).
\end{thm*}

Influenced by the works \cite{CY,LY,Y} of Yau’s school and the work of Gidas-Spruck \cite{GS}, in \cite{LJY}, J.Y.Li first derived Yau’s gradient estimate and local \({L}^{\infty }\) estimate for equation \eqref{LE} on the manifolds with Ricci curvature bounded below for \(1 \leq  \alpha  < \; \frac{n}{n - 2}\) and \(n \geq  4\) . Notice that the gradient estimate is strictly stronger than the Liouville theorem for equation \eqref{LE} and hence he provided a new proof of Gidas-Spruck's Liouville theorem when \(1 \leq  \a < \frac{n}{n - 2}\) and \(n \geq  4\) . Later, Y. Y. Yang \cite{Yang} obtained the same estimate for \eqref{LE} with \(\alpha<0\).
 Recently, Peng-Wang-Wei \cite{PWW}, Wang-Wei \cite{WW} and Lu \cite{LU1} have successively refined the exponent for which Yau's gradient estimate holds until reaching \(\alpha  < \frac{n + 3}{n - 1}\) . Using the Keller--Osserman comparison method, Bidaut-Véron--García-Huidobro--Véron \cite{BGV} obtained a gradient estimate for \eqref{LE} on Euclidean spaces with \(\alpha  \in  \left\lbrack  {1,\frac{n + 3}{n - 1}}\right)\) . Using Moser’s iteration method from Wang-Zhang \cite{WZ}, He-Wang-Wei \cite{HWW} also obtain the same results and their $p$-Laplacian generalization. At present, for the subcritical equation \eqref{LE} on manifolds with Ricci curvature bounded below, Yau's gradient estimates and the uniform bound estimates of solutions have not yet been fully established.

Inspired by Gidas-Spruck's Liouville theorem, one would naturally conjecture that the stronger Yau's gradient estimate and universal bounds of solutions might also hold for all subcritical equations of \eqref{LE}. We formulate this long-standing question as follows.

	\begin{problem}\label{ques:1}
		Do Cheng--Yau's gradient estimate and the universal bounds estimate for solutions to equation \eqref{LE} hold for all subcritical cases?
	\end{problem}
	
Concretely, for any real number $n\in[1,\infty)$, we define
\begin{eqnarray}
	p_S(n)=\begin{cases}
		\infty    \quad\qquad\qquad\text{if}\quad n\in[1,2]\nonumber\\
		\frac{n+2}{n-2}\qquad\qquad\,\,\,\text{if}\quad n\in(2,\infty).\nonumber
	\end{cases}
\end{eqnarray}

The main result of this paper is to solve the above problem.
	\begin{theorem}\label{mthm}
		Let $(\M^n,g)$ be an $n$-dimensional complete Riemannian manifold with $Ric\ge -Kg$ and $K\ge 0$.  If $u$ is a positive  solution of \eqref{LE} on $B(x_0,2R)\subset\M^n$ with $\a\in(-\infty,p_S(n))$, then 
		\begin{equation}\label{strong1}
			\sup\limits_{B(x_0,R)}\left(\frac{|\nabla u|^2}{u^2}+u^{\alpha-1}\right)\le C(n,\a)\left(K+\frac{1}{R^2}\right).
		\end{equation}	
		
	\end{theorem}

	From Theorem \ref{mthm}, we immediately obtain Gidas--Spruck's Liouville theorem, universal $L^{\infty}$ estimates for $\alpha\in (1,p_S(n))$ (when $\alpha<1$, it implies a universal lower bound estimate), Yau's gradient estimate and its corollary---Harnack inequality. Even in the Ricci non-negative case, estimate \eqref{strong1} is new. Simultaneously, estimate \eqref{strong1} is also sharp in Euclidean case by singular exterior solutions, see \cite{PQS,QS,SZ}.

	\begin{remark}\label{r11}
		We mention that the positive constant $C(n,\alpha)$  can be calculated. When  $\alpha$  is in some left neighborhood of  $\frac{n+2}{n-2}$ with $n\ge 3$,  $	C(n,\alpha)$  has the following expression:
		$$
		C(n,\alpha)=\gamma(n)\left(\frac{n+2}{n-2}-\alpha\right)^{-\theta(n)},
		$$
		where $\gamma(n),\theta(n)$ are positive constants depending only on  $n$.
	\end{remark}
	
	As stated earlier, the Harnack inequality is a natural corollary.
		\begin{corollary}\label{HK}
		Let $(\M^n,g)$ be an $n$-dimensional complete Riemannian manifold with $Ric\ge -Kg$ and $K\ge 0$.  If $u$ is a positive  solution of \eqref{LE} on $B(x_0,2R)\subset\M^n$ with $\a\in(-\infty,p_S(n))$, then 
		\begin{equation}\label{shk}
			\sup\limits_{B(x_0,R)}u\le Ce^{C\sqrt{K}R}\inf\limits_{B(x_0,R)}u,
		\end{equation}	
		where $C=C(n,\alpha)$.
	\end{corollary}

		Before we further elaborate on the contributions of this paper, we provide a historical remark on the research of equation \eqref{LE}.
		
		\begin{remark}
			 When considering the same equation \eqref{LE} in Euclidean space, in addition to the original integral estimation method by Gidas-Spruck, the moving plane method and moving sphere method also provide elegant proofs for the Liouville theorem of Gidas-Spruck. These methods have also successfully classified the solutions to the critical equations (which is a special Yamabe equation, see the classical works of Schoen \cite{Sch} and Lee-Parker \cite{LP}), and readers may refer to the classic works of Aleksandrov \cite{A}, Serrin \cite{Se}, Gidas-Ni-Nirenberg \cite{GNN}, Caffarelli-Gidas-Spruck \cite{CGS}, Chen-Li \cite{CL} and Li-Zhang \cite{LZ} and their subsequent studies. From Liouville's theorem to obtaining uniform upper bounds and gradient estimates for local solutions in Euclidean spaces, there are currently two methods applicable to all subcritical exponents: One approach is to construct a blow-up sequence and then use Liouville’s theorem to derive a contradiction, thereby obtaining the \({L}^{\infty }\) estimate of the solution. This method usually requires the exponent \(\alpha  > 1\) . The other classical method is to first obtain the Harnack inequality for local solutions in the subcritical case based on the important identity of Gidas and Spruck, and then derive the upper bound estimate of the solution using standard elliptic inequalities. These methods can be found in the important works of Gidas-Spruck \cite{GS0}, Dancer \cite{D}, Serrin-Zou \cite{SZ} and Poláčik-Quittner-Souplet \cite{PQS}. Both the Liouville theorems on Euclidean spaces and the methods for local estimates have been included in Quittner-Souplet's excellent book \cite{QS}. We recommend that readers refer to it for more details.
		\end{remark}

	\subsection{New contribution of the paper}
The main contribution of this paper in studying Problem \ref{ques:1} lies in the discovery of  new transformations and the corresponding new auxiliary functions. Although we also use the maximum principle like many previous works, our argumentation method is completely new. Below, we  will specifically explain the difficulty of the problem and our strategy.

Previous works usually use the square of the norm of the logarithmic gradient or some of its variants as auxiliary functions, such as
$$\frac{|\nabla u|^2}{u^2}\quad \text{or}\quad \frac{|\nabla u|^2}{u^2}+du^{\a-1},$$
where $u$ is  a solution of \eqref{LE} and $d>0$ is undetermined. Meanwhile, exponential transformations $w=u^{-\b}$  can often better expand the range of indices (these transformations have been widely used in previous works of Obata \cite{O}, Gidas--Spruck \cite{GS}, J.Y.Li \cite{LJY} and their subsequent studies), but they cannot cover all subcritical indices. 
Inspired by \cite{GS}, the auxiliary function ($\g,d$ are undetermined)
$$u^{\g}\left(\frac{|\nabla u|^2}{u^2}+du^{\a-1}\right)$$
seems to be a reasonable candidate because it plays an important role in the derivation of the Gidas--Spruck integral estimate and the Liouville theorem.
However, if we directly use this function as an auxiliary function, we cannot obtain Yau's logarithmic gradient estimate. The specific difficulties are reflected in the following two technical aspects:\\
$\bullet$
After applying the Bochner formula, the condition of a lower bound on curvature becomes a significant obstacle to the argument, unless it is strengthened to non-negative curvature;\\
$\bullet$ To derive the logarithmic gradient estimate from the \( L^\infty \) estimate of the auxiliary function $$u^{\g}\left(\frac{|\nabla u|^2}{u^2}+du^{\a-1}\right),$$ one actually needs the solution to satisfy the Harnack inequality a priori. However, up to now, the  Harnack inequality \eqref{shk} for the subcritical equation \eqref{LE} on manifolds has been lacking.

To overcome these difficulties, we introduce the following auxiliary function with parameter \( \varepsilon \) and transformation (for more detailed formulas, see Section 2):
$$
(u+\e)^{-\b\g}\left(\frac{|\nabla w|^2}{w^2}+du^{\a-1}\right) \quad\text{with} \quad w=u^{-\b};
$$
and 
$$
w^{\g}\left(\frac{|\nabla w|^2}{w^{2}}+du^{\a-1}\right)\quad \text{with} \quad w=(u+\e)^{-\b}.
$$
By computing their auxiliary functions, we have successfully derived their \( L^\infty \) estimates. Of particular importance is that the constants in these \( L^\infty \) estimates are independent of \( \varepsilon \).
Finally, through a very delicate choice of \( \varepsilon \) (\( \varepsilon \) is a root of an algebraic equation related to the variables \( n \), \( \a \), \( K \) and \( R \); for details, see \eqref{47}), we derived the desired estimate \eqref{strong1} from the \( L^\infty \) estimate of the above auxiliary functions with \( \varepsilon \).

	\subsection{Notation and the organization of the paper}\label{notation-sec}

	We use $(\M^n,g)$ to denote an $n$-dimensional complete Riemannian manifold with Riemannian metric $g$ and $B(x_0,R)\subset\M^n$ to denote the open ball in $\M^n$ with center $x_0$ and radius $R$. 
	Let $\Omega\subset\M^n$ be an open subset, we define $ C^k(\Omega) $ to be the space of all real-valued functions on $ \Omega $ that are $ k $ times continuously differentiable. If $ u \in C^1(\Omega) $, we use $\nabla u$ to denote its gradient; and if $ u \in C^2(\Omega) $, we use $\nabla^2 u$ to denote its Hessian tensor. Throughout the entire paper, the solutions of the equations we mentioned all refer to classical solutions.
	We use $|T|$ to denote the norm of tensor $T$ under the Riemannian metric $g$.
 We use $X\le C(\cdot)Y$ to denote the estimate $X \leq CY$ for some positive constant $C$ which depends only on the parameters listed in the parentheses.
	For a set \( E \), we use \( \chi_E \) to denote its characteristic function, that is,
	\begin{eqnarray}
	\chi_E(x)=\begin{cases}
			1    \quad\qquad\qquad\text{if}\quad x\in E\nonumber\\
			0\qquad\qquad\quad\text{if}\quad x\notin E.\nonumber
		\end{cases}
	\end{eqnarray}

	The rest of this paper is organized as follows. In Section \ref{s2}, we introduce two key auxiliary functions and provide the elliptic equations associated with these functions, which serve as the foundations for the subsequent arguments.  In Section \ref{s3},  
	we have established several lemmas that will play a crucial role in the proof of Theorem \ref{mthm} later.
Then	Section \ref{s4} will be devoted to the proof of Theorem \ref{mthm}.  Section \ref{s5} generalizes Theorem \ref{mthm} to an equation with gradient term on manifolds with Bakry-Émery Ricci curvature bounded below.
	Simultaneously, in Section \ref{s5}, an example is given to illustrate that our estimate is optimal.

After the initial preprint of the present paper \cite{LU0} was posted to the arXiv (the results of the revised manuscript are fully contained in \cite{LU0}), Wu \cite{W} gave a simplified proof of the Gidas–Spruck Liouville theorem by setting $\varepsilon=0$ in \eqref{F} and using lower bound estimates for superharmonic functions, with a focus on Liouville-type theorems (see also Lu \cite{LU2}). Within half a year after completing this work at USTC, we generalized our results to more general equations on metric measure spaces in \cite{LU3}.

	\subsection{Acknowledgments}
		We would like to thank Professor Jiayu Li and Professor Youde Wang for many helpful and insightful discussions on this topic. We also warmly thank the referee for carefully reading our manuscript and for the helpful suggestions, which have improved the paper.
	The author is supported by the Scientific Research Startup Project of Jiangsu Normal University (Project No. 24XFRS051), the General Program of Basic Scientific Research in Institutions of Higher Education of Jiangsu Province (Grant No. 25KJB110002) and the National Natural Science Foundation of China (Grant No. 12526552).

	\section{Elliptic equations of auxiliary functions} \label{s2}

	The core of the section is the construction of two different auxiliary functions related to solutions and the computation of their elliptic equations respectively. 
	Let $u$ be a positive solution of \eqref{LE} on $B(x_0,2R)\subset\M^n$. For $\b\neq 0$, we set
	\begin{equation}\label{tran}
		w=u^{-\beta},
	\end{equation}
	then the chain rule gives the equivalent equation of \eqref{LE}:
	\begin{equation}\label{wle}
		\Delta w=\left(1+\frac{1}{\beta}\right)\frac{|\nabla w|^2}{w}+\b wu^{\a-1}.
	\end{equation}
	For undetermined real numbers $\g\ge 0,\e\ge 0,d\ge 0$, we define \emph{the first type} auxiliary function:
	\begin{equation}\label{F}
		F=(u+\e)^{-\b\g}\left(\frac{|\nabla w|^2}{w^2}+du^{\a-1}\right).
	\end{equation}
	For an undetermined $\e>0$, we also need another transformation:
	\begin{equation}\label{tran2}
		w=(u+\e)^{-\b}.
	\end{equation}
 Under transformation \eqref{tran2}, equation \eqref{LE} becomes
	\begin{equation}\label{wlee}
		\d w=\left(1+\frac{1}{\b}\right)\frac{|\nabla w|^2}{w}+\b wu^{\a-1}\times\frac{u}{u+\e}.
	\end{equation}
	We define \emph{the second type} auxiliary function:
	\begin{equation}\label{G}
		G=w^{\g}\left(\frac{|\nabla w|^2}{w^{2}}+du^{\a-1}\right),
	\end{equation}
	where $\g, d\ge0$ are undetermined real numbers.

	\begin{remark}\label{r21}
		The first type and the second type of auxiliary functions work for dimensions greater than or equal to 4 and less than 4 respectively. This is determined when we make the lower bound estimate for the coefficient of the principal term in our proof, because when we truncate the solution, we need a positive lower bound of the principal term coefficient that is independent of $\e$. This forces us to use the second type of auxiliary function instead of the first type when the dimension is less than 4. The technical details can be found in Subsection \ref{s33}.
	\end{remark}
	
	The following two lemmas give elliptic equations about auxiliary functions $F,G$.
	\begin{lemma}\label{k1}
		Let $u$ be a positive solution of \eqref{LE} on $B(x_0,2R)\subset\M^n$, $w$ be defined by \eqref{tran} with $\b\neq 0$ and $F$ be its first type auxiliary function as \eqref{F} with $\g,\e,d\ge 0$. Then we have
		\begin{eqnarray}
			(u+\e)^{\b\g}\Delta F&=&2w^{-2}\left|\n2-\frac{\d w}{n}g\right|^2+2w^{-2}Ric(\nabla w,\nabla w)\nonumber\\
			&&+2\left(\frac{1}{\b}-1+\g\frac{u}{u+\e}\right)(u+\e)^{\b\g}\left\langle\nabla F,\nabla \ln w\right\rangle\nonumber\\
			&&+U\frac{|\nabla w|^4}{w^4}+V\frac{|\nabla w|^2}{w^2}u^{\a-1}+Wu^{2(\a-1)}
		\end{eqnarray}
		on $B(x_0,2R)$, where
		\begin{displaymath}
			\begin{aligned}
				U=&\,\frac{2}{n}\left(1+\frac{1}{\b}\right)^2+\left(\frac{\g}{\b}-\g^2\right)\frac{u^2}{(u+\e)^2}+2\left(1-\frac{1}{\b}\right)\frac{\g u}{u+\e}-2,\\
				V=&\,\frac{4}{n}(1+\b)+2\left(1-\a\right)+d\left(\frac{\a(\a-1)}{\b^2}-\frac{2(\a-1)}{\b}\right)\nonumber\\
				&+\frac{\g u}{u+\e}\left(\b+d\left(\left(\frac{1}{\b}-\g\right)\frac{u}{u+\e}+2-\frac{2}{\b}\right)\right),\\
				W=&\,\frac{2\b^2}{n}+d\left(\frac{\b\g u}{u+\e}+1-\a\right).\nonumber
			\end{aligned}
		\end{displaymath}

	\end{lemma}
	\begin{proof}
		By calculus rule, we have
		\begin{eqnarray}\label{11}
			\d F&=&\d((u+\e)^{-\b\g})\times\left(\frac{|\nabla w|^2}{w^2}+du^{\a-1}\right)\nonumber\\
			&&+2\left\langle\nabla (u+\e)^{-\b\g},\nabla\left(\frac{|\nabla w|^2}{w^2}\right)\right\rangle\nonumber\\
			&&+2d\left\langle\nabla (u+\e)^{-\b\g},\nabla\left(u^{\a-1}\right)\right\rangle\nonumber\\
			&&+(u+\e)^{-\b\g}\left(\d\left(\frac{|\nabla w|^2}{w^2}\right)+d\d\left(u^{\a-1}\right)\right).
		\end{eqnarray}
		Now, we concretely compute the terms in \eqref{11} as follows.
		\begin{equation}\label{12}
			\d((u+\e)^{-\b\g})=\b\g (u+\e)^{-\b\g-1}u^{\a}+\b\g(\b\g+1)(u+\e)^{-\b\g-2}|\nabla u|^2,
		\end{equation}
		\begin{eqnarray}\label{13}
			\d u^{\a-1}=(\a-1)u^{\a-2}+(\a-1)(\a-2)u^{a-3}|\nabla u|^2,
		\end{eqnarray}
		\begin{eqnarray}\label{14}
			\left\langle\nabla (u+\e)^{-\b\g},\nabla(u^{\a-1})\right\rangle=-\b\g(\a-1) (u+\e)^{-\b\g-1}u^{\a-2}|\nabla u|^2,
		\end{eqnarray}
		\begin{eqnarray}\label{15}
			\left\langle\nabla (u+\e)^{-\b\g},\nabla\left(\frac{|\nabla w|^2}{w^2}\right)\right\rangle&=&\frac{\g u}{(u+\e)^{\b\g+1}}\left(-2\frac{|\nabla w|^4}{w^4}+\frac{\left\langle\nabla w,\nabla|\nabla w|^2\right\rangle}{w^3}\right).\nonumber\\
			&&
		\end{eqnarray}
		Then we deal with the term $\d\left(\frac{|\nabla w|^2}{w^2}\right)$.
		\begin{eqnarray}\label{16}
			\d\left(\frac{|\nabla w|^2}{w^2}\right)&=&\d(w^{-2})|\nabla w|^2+w^{-2}\d|\nabla w|^2+2\left\langle\nabla|\nabla w|^2,\nabla(w^{-2})\right\rangle\nonumber\\
			&=&\left(-2w^{-3}\d w+6w^{-4}|\nabla w|^2\right)|\nabla w|^2\nonumber\\
			&&-4w^{-3}\left\langle\nabla w,\nabla|\nabla w|^2\right\rangle+w^{-2}\d|\nabla w|^2\nonumber\\
			&=&-2w^{-3}|\nabla w|^2\left(\left(1+\frac{1}{\b}\right)\frac{|\nabla w|^2}{w}+\b wu^{\a-1}\right)\nonumber\\
			&&+6w^{-4}|\nabla w|^4-4w^{-3}\left\langle\nabla w,\nabla|\nabla w|^2\right\rangle+w^{-2}\d|\nabla w|^2\nonumber\\
			&=&\left(4-\frac{2}{\b}\right)w^{-4}|\nabla w|^4-2\b w^{-2}|\nabla w|^2u^{\a-1}\nonumber\\
			&&-4w^{-3}\left\langle\nabla w,\nabla|\nabla w|^2\right\rangle+2w^{-2}|\nabla^2 w|^2+2w^{-2}Ric(\nabla w,\nabla w)\nonumber\\
			&&+2w^{-2}\left\langle\nabla w,\nabla\left(\left(1+\frac{1}{\b}\right)\frac{|\nabla w|^2}{w}+\b wu^{\a-1}\right)\right\rangle\nonumber\\
			&=&2\left(\frac{1}{\b}-1\right)\frac{\left\langle\nabla w,\nabla|\nabla w|^2\right\rangle}{w^{3}}+\left(2-\frac{4}{\b}\right)w^{-4}|\nabla w|^4\nonumber\\
			&&+2\left(1-\a\right)u^{\a-1}\frac{|\nabla w|^2}{w^2}+2\frac{|\nabla^2 w|^2}{w^2}+2Ric(\nabla \ln w,\nabla\ln w)\nonumber\\
			&&
		\end{eqnarray}
		on $B(x_0,2R)$,
		where we use equation \eqref{wle} and the chain rule for the second equation and the third equation, Bochner formula and equation \eqref{wle} for the fourth equation and last equation respectively.

		Through direct calculations, we have the following two identities.
		
		\begin{eqnarray}\label{19}
			\frac{\left\langle\nabla w,\nabla|\nabla w|^2\right\rangle}{w^3}&=&(u+\e)^{\b\g}\left\langle\nabla F,\nabla\ln w\right\rangle+\left(2-\frac{\g u}{u+\e}\right)\frac{|\nabla w|^4}{w^4}\nonumber\\
			&&-d\left(\frac{\g u}{u+\e}+\frac{1}{\b}\left(1-\a\right)\right)\frac{|\nabla w|^2}{w^2}u^{\a-1},
		\end{eqnarray}
		and
		\begin{eqnarray}\label{id}
			|\nabla^2 w|^2&=&\left|\n2-\frac{\d w}{n}g\right|^2+\frac{1}{n}(\d w)^2\nonumber\\
			&=&\left|\n2-\frac{\d w}{n}g\right|^2+\frac{1}{n}\left(\left(1+\frac{1}{\beta}\right)\frac{|\nabla w|^2}{w}+\b wu^{\a-1}\right)^2,
		\end{eqnarray}
	Substitute \eqref{12}-\eqref{16} into \eqref{11} first, then we use  \eqref{19} and \eqref{id} to substitute the terms $w^{-3}\left\langle\nabla w,\nabla|\nabla w|^2\right\rangle$ and $|\nabla^2 w|^2$.	Finally, we combine the coefficients of the principal terms and finish the proof.

	\end{proof}

Following a calculation process similar to the proof of Lemma \ref{k1}, it is not difficult for us to derive the elliptic equation of the second type auxiliary function \( G \). We directly present its equation while omitting its proof.
\vspace{2mm}

	\begin{lemma}\label{k2}
		Let $u$ be a positive solution of \eqref{LE} on $B(x_0,2R)\subset\M^n$, $w$ be defined as \eqref{tran2} with $\b\neq 0$ and $G$ be its second type auxiliary function as \eqref{G} with $\e,d>0$ and $\g=1$. Then we have
		\begin{eqnarray}\label{k22}
			w^{-1}\Delta G&=&2w^{-2}\left|\n2-\frac{\d w}{n}g\right|^2\nonumber\\
			&&+2w^{-2}Ric(\nabla w,\nabla w)+\frac{2}{\b}w^{-\g}\left\langle\nabla G,\nabla \ln w\right\rangle\nonumber\\
			&&+U\frac{|\nabla w|^4}{w^4}+V\frac{|\nabla w|^2}{w^2}u^{\a-1}+Wu^{2(\a-1)}
		\end{eqnarray}
		on $B(x_0,2R)$, where
		\begin{displaymath}
			\begin{aligned}
				U=&\,\left(\frac{2}{n}\left(1+\frac{1}{\b}\right)-1\right)\left(1+\frac{1}{\b}\right),\\
				V=&\,\left(\frac{4}{n}(1+\b)+2+\b\right)\frac{u}{u+\e}-2\a+\frac{2d}{\b}\left(\frac{\a-1}{\b}\frac{u+\e}{u}-1\right)\\
				&+d\left(\frac{(\a-1)(\a-2)}{\b^2}\frac{(u+\e)^2}{u^2}+1+\frac{1}{\b}-\frac{2(\a-1)}{\b}\frac{u+\e}{u}\right),\\
				W=&\,\frac{2\b^2}{n}\left(\frac{u}{u+\e}\right)^2+d\left(\frac{\b\g u}{u+\e}+1-\a\right).\nonumber
			\end{aligned}
		\end{displaymath}

	\end{lemma}

	\section{Some preliminary lemmas}\label{s3}
	In this section, we collect lemmas used in the proof of Theorem \ref{mthm}. First,
	as in Yau \cite{Y}, Cheng-Yau \cite{CY} or Li-Yau \cite{LY},  by the Laplacian comparison theorem, we have a good cut-off function.
	\begin{lemma}\label{cut off}
		Let $(\M^n,g)$ be an $n$-dimensional complete Riemannian manifold with $Ric\ge-Kg$, where $K\ge 0$. Then for any $R>0$, we have a Lipschitz function $\Phi$ on $B(x_0,2R)$ such that\\
		\rm(i) $\Phi(x)=\phi(d(x_0,x))$, where  $d(x_0,\cdot)$ is the distance function from $x_0$ and  $\phi$ is a non-increasing function on $[0,\infty)$ and
		\begin{eqnarray}
			\Phi(x)=
			\begin{cases}
				1\qquad\qquad\qquad\,\text{if}\qquad x\in B(x_0,R)\nonumber\\
				0\qquad\qquad\qquad\,\text{if}\qquad x\in B(x_0,2R)\setminus B(x_0,\frac{3}{2}R).
			\end{cases}
		\end{eqnarray}
		(ii) On $\{x\in B(x_0,2R):\Phi(x)>0\}$,
		\begin{equation}
			\frac{|\nabla \Phi|}{\Phi^{\frac{1}{2}}}\le\frac{C}{R}.\nonumber
		\end{equation}
		(iii)
		\begin{equation}
			\Delta \Phi\ge-\frac{C\sqrt{nK}}{R} \coth\left( \sqrt{\frac{K}{n}}R\right)-\frac{C}{R^2}\ge- C(n)\frac{\sqrt{K}R+1}{R^2} \nonumber
		\end{equation}
		holds	on $B(x_0,2R)$ in the distribution sense and pointwise outside the cut locus of $x_0$. Here, $C$ is a universal constant.
	\end{lemma}

	\subsection{A partial result}
	
		By directly taking the squared norm of the logarithmic gradient as the auxiliary function, we can obtain a partial result of Theorem \ref{mthm}.

	\begin{lemma}\label{ps}
		Let $(\M^n,g)$ be an $n$-dimensional complete Riemannian manifold with $Ric\ge -Kg$ and $K\ge 0$. If $u$ is a positive solution of \eqref{LE} on $B(x_0,2R)\subset \M^n$ with $\a<\frac{n+3}{n-1}$ and $n\ge2$ or $\a\in\R$ and $n=1$, then 
		\begin{equation}
			\sup\limits_{B(x_0,R)}\frac{|\nabla u|^2}{u^2}\le C(n,\a)\left(K+\frac{1}{R^2}\right),
		\end{equation}

	\end{lemma}
	
	\begin{proof}
	Define $w$ by \eqref{tran} with $\b=\frac{2}{n-1}$ and $n\ge 2$ or $\b=|\a|+1$ and $n=1$,  and let $F$ be its first type auxiliary function as \eqref{F} with $\g,\e,d=0$. Then we have ($n\ge 2$)
	\begin{eqnarray}\label{32}
		\Delta F&=&2w^{-2}\left|\n2-\frac{\d w}{n}g\right|^2+2w^{-2}Ric(\nabla w,\nabla w)\nonumber\\
		&&+2\left(\frac{1}{\b}-1\right)\left\langle\nabla F,\nabla \ln w\right\rangle+\frac{2}{n\b^2}\frac{|\nabla w|^4}{w^4}+\frac{2\b^2}{n}u^{2(\a-1)}\nonumber\\
		&&+\left(\frac{4}{n}(1+\b)+2(1-\a)\right)\frac{|\nabla w|^2}{w^2}u^{\a-1}
	\end{eqnarray}
	on $B(x_0,2R)$. By the basic inequality, for $\delta\in(0,2)$, we have
	\begin{equation}\label{33}
		\frac{2-\delta}{n\b^2}\frac{|\nabla w|^4}{w^4}+\frac{(2-\delta)\b^2}{n}u^{2(\a-1)}\ge \frac{4-2\delta}{n}\frac{|\nabla w|^2}{w^2}u^{\a-1}.
	\end{equation}
	Since \( \a < \frac{n+3}{n-1}\),  there exists a \( \delta=\delta(n,\a)>0 \) such that
		\begin{equation}\label{34}
		\frac{4}{n}(1+\b)+2(1-\a)+\frac{4-2\delta}{n}\ge 0.
	\end{equation}
	Combining the \eqref{32}, \eqref{33} and \eqref{34} with the curvature condition of the manifold, we obtain
		\begin{eqnarray}\label{35}
		\Delta F&\ge&\frac{\delta}{n\b^2}\frac{|\nabla w|^4}{w^4}+\frac{\delta\b^2}{n}u^{2(\a-1)}\nonumber\\
		&&-2K\frac{|\nabla w|^2}{w^2}+2\left(\frac{1}{\b}-1\right)\left\langle\nabla F,\nabla \ln w\right\rangle
	\end{eqnarray}
	on $B(x_0,2R)$. From inequality \eqref{35}, it is not difficult for us to obtain the desired estimation, with specific details as follows.

	Let $A:=\Phi F$ be the auxiliary function on $B(x_0,2R)$, where $\Phi$ is a cut-off function as in Lemma \ref{cut off}. Without loss of generality, we may assume that $A$ attains its positive maximum value in $B(x_0,\frac{3}{2}R)$ or the proof is finished. We also assume that the maximal value point $x$ is outside of the cut locus of $x_0$ or one uses Calabi's argument as Yau \cite{Y}, Cheng-Yau \cite{CY} or Li-Yau \cite{LY}.
	By chain rule, at point $x$, we can see
	\begin{eqnarray}\label{aeq}
		0\ge \d A&=&F\d \Phi+\Phi\d F+2\left\langle\nabla \Phi,\nabla F\right\rangle\nonumber\\
		&=&F\d \Phi+\Phi\d F-2A\frac{|\nabla \Phi|^2}{\Phi^2}.
	\end{eqnarray}
	By inequality \eqref{35}, \eqref{aeq} and chain rule  again, we have
	\begin{eqnarray}\label{42}
		0&\ge& \frac{\delta}{n\b^2}\Phi F^2-2KA+F\Delta \Phi-2A\frac{|\nabla \Phi|^2}{\Phi^2}\nonumber\\
		&&-2\left(\frac{1}{\b}-1\right)F\left\langle\nabla \Phi,\nabla \ln w\right\rangle.
	\end{eqnarray}
Define $L=\frac{\delta}{n\b^2}$, then the Cauchy-Schwarz inequality gives
	\begin{eqnarray}\label{bb}
		-2\left(\frac{1}{\b}-1\right)F\left\langle\nabla \Phi,\nabla \ln w\right\rangle\ge-\left(\frac{1}{\b}-1\right)^2\frac{2}{L}\frac{|\nabla\Phi|^2}{\Phi}F-\frac{L}{2}\Phi F^2.
	\end{eqnarray}
	Substituting \eqref{bb} into \eqref{42} yields
	\begin{eqnarray}\label{4300}
		0&\ge& \frac{L}{2}\Phi F^2+ F\Delta \Phi-2KA-2\left(\frac{1}{L}\left(\frac{1}{\b}-1\right)^2+1\right)\frac{|\nabla \Phi|^2}{\Phi^2}A.
	\end{eqnarray}
	Multiplying $\Phi(x)$ on both sides of \eqref{4300} and using the property of $\Phi$ in Lemma \ref{cut off}, we have
	\begin{equation}\label{45ps}
		A(x)\le C(n,\a)\left(K+\frac{1}{R^2}\right),
	\end{equation}
	which implies
	\begin{equation}
		\sup\limits_{B(x_0,R)}\frac{|\nabla u|^2}{u^2}\le C(n,\a)\left(K+\frac{1}{R^2}\right).
	\end{equation}
For the case of \( n = 1 \), we can directly obtain an inequality similar to inequality \eqref{35} by Lemma \ref{k1}. The subsequent part of the argument is almost identical, and thus we have completed the proof.

	\end{proof}

	\subsection{Control lemmas}
	
When $\a>1$, we do not have gradient estimates. Instead, we have  weaker control lemmas.
	\begin{lemma}\label{cl1}
			Let $(\M^n,g)$ be an $n$-dimensional complete Riemannian manifold with $Ric\ge -Kg$ and $K\ge 0$.  If $u$ is a positive  solution of \eqref{LE} on $B(x_0,2R)\subset\M^n$ with $\a\in\R$, then 
		\begin{equation}\label{312}
			\sup\limits_{B(x_0,R)}\frac{|\nabla u|^2}{u^2}\le C(n,\a)\left(K+\frac{1}{R^2}+\sup\limits_{B(x_0,2R)}u^{\alpha-1}\right).
		\end{equation}	
		Furthermore, $C(n,\a)=C(n)\a$ when $\a>1$.
	\end{lemma}
	\begin{proof}
	According to Lemma \ref{ps}, we only need to prove the case where \( \a > 1 \). For this, 	we define $w$ by \eqref{tran} with $\b=\frac{1}{n}$  and let $F$ be its first type auxiliary function as \eqref{F} with $\g,\e,d=0$. Then we have
	\begin{eqnarray}\label{321}
		\Delta F&\ge&2n\frac{|\nabla w|^4}{w^4}+\frac{2}{n^3}u^{2(\a-1)}+2\left(n-1\right)\left\langle\nabla F,\nabla \ln w\right\rangle\nonumber\\
		&&-2\a\frac{|\nabla w|^2}{w^2}u^{\a-1}-2K\frac{|\nabla w|^2}{w^2}
	\end{eqnarray}
	on $B(x_0,2R)$.  Starting from inequality \eqref{321}, we can obtain estimate \eqref{312} through arguments similar to those in the latter part of Lemma \ref{ps} and we omit these redundant details.

	\end{proof}

	\begin{lemma}\label{cl2}
		Let $(\M^n,g)$ be an $n$-dimensional complete Riemannian manifold with $Ric\ge -Kg$ and $K\ge 0$.  If $u$ is a positive  solution of \eqref{LE} on $B(x_0,2R)\subset\M^n$ with $\a\in\R$, then 
		\begin{equation}\label{314}
			\sup\limits_{B(x_0,R)}u^{\alpha-1}\le C(n,\a)\left(K+\frac{1}{R^2}+\sup\limits_{B(x_0,2R)}\frac{|\nabla u|^2}{u^2}\right).
		\end{equation}	
			Furthermore, $C(n,\a)=C(n)\a^2$ when $\a>1$.
	\end{lemma}
	\begin{proof}
		We define $w$ by \eqref{tran} with $\b=n(|\a|+1)$  and let $F$ be its first type auxiliary function as \eqref{F} with $\g,\e=0,d=1$. Then we have
		\begin{eqnarray}\label{322}
			\Delta F&\ge&-2\frac{|\nabla w|^4}{w^4}+2n(\a^2+1)u^{2(\a-1)}+2\left(\frac{1}{\b}-1\right)\left\langle\nabla F,\nabla \ln w\right\rangle\nonumber\\
			&&+(2+2|\a|)\frac{|\nabla w|^2}{w^2}u^{\a-1}-2K\frac{|\nabla w|^2}{w^2}
		\end{eqnarray}
		on $B(x_0,2R)$.  Starting from inequality \eqref{322}, we can obtain estimate \eqref{314} through arguments similar to those in the latter part of Lemma \ref{ps} and we omit these redundant details.
		
	\end{proof}

	\subsection{Estimation of the leading coefficients}\label{s33}
	As mentioned in Remark \ref{r21}, the lower bound estimation of the leading term coefficient constitutes the core part of the proof of Theorem \ref{mthm}---specifically, the following two lemmas. The key point lies in the independence of the lower bound of the coefficient of the fourth power of the gradient norm from \( \e \).

	\begin{lemma}\label{ab}
		Let $n\ge 4$, $\e>0, \g=1$ and $\a\in[\frac{n+3}{n-1},p_S(n))$. The functions $U, V, W$ are defined as in Lemma \ref{k1}. Then there exist constants $\b=\b(n,\a)\in(0,\frac{2}{n-2})$, $d=d(n,\a)>0$,  $L=L(n,\a)>0$ and $M=M(n,\a)>0$ such that for any $\e>0$,
		\begin{equation}
			U\ge U_0>0,\nonumber
		\end{equation}
		\begin{equation}
			V\ge V_0-M\chi_{\{x\in B(x_0,2R): u(x)< L\e\}},\nonumber
		\end{equation}
		\begin{equation}
			W\ge W_0-M\chi_{\{x\in B(x_0,2R): u(x)< L\e\}},\nonumber
		\end{equation}
		where $U_0,V_0,W_0$ are positive constants which only depend on $n, \a$.	
	\end{lemma}
	\begin{proof}
		Let $\g=1$ and $U, V, W$ be defined as in Lemma \ref{k1}, then	
		\begin{eqnarray}
			U&=&\left(\frac{1}{\b}-1\right)\left(\frac{u}{u+\e}-2\right)\frac{u}{u+\e}+\frac{2}{n}\left(1+\frac{1}{\b}\right)^2-2\nonumber\\
			&\ge&\left(\frac{2}{n}\left(1+\frac{1}{\b}\right)-1\right)\left(1+\frac{1}{\b}\right)>0,\nonumber
		\end{eqnarray}
		when $\frac{1}{\b}-1\ge0$ and $\b<\frac{2}{n-2}$, which is equivalent to $\b<\frac{2}{n-2}$ if $n\ge 4$. This is the only place where we use the condition that \( n \geq 4 \).
		
		We define
		\begin{equation}\label{uc}
			U_c=\left(\frac{2}{n}\left(1+\frac{1}{\b}\right)-1\right)\left(1+\frac{1}{\b}\right).
		\end{equation}
		On the other hand, from the formulas of $V$ and $W$, we know that
		\begin{eqnarray}
			V&=&V_c+d\left(\frac{1}{\b}-1\right)\frac{\e^2}{(u+\e)^2}-\frac{\b\e}{u+\e},\nonumber\\
			W&=&W_c-\frac{d\b\e}{u+\e}\nonumber,
		\end{eqnarray}
		where 
		\begin{eqnarray}
			V_c&=&\frac{4}{n}(1+\b)+2\left(1-\a\right)+\b+\frac{d}{\b^2}\left(\a-\b-1\right)\left(\a-\b\right),\label{vc}\\
			W_c&=&\frac{2}{n}\b^2-d\left(\a-\b-1\right)\label{wc}.
		\end{eqnarray}
		
	For a to-be-determined \(\b \in (0, \frac{2}{n-2}) \),	we set $d_0=\frac{2\b^2}{n(\a-1-\b)}>0$, because according to the conditions, we have $\a\ge\frac{n+3}{n-1}>\frac{n}{n-2}>1+\b$ if $n\ge 4$.  If we choose $d=d_0$ first, direct algebraic computation shows that:
	
	 $$\text{when}\quad \b\in\left(2\left(\frac{n-1}{n+2}\a-1\right),\a-\frac{n-1}{n-2}\right),\qquad V_0>0.$$ Then we choose a fixed $\b\in\left(2\left(\frac{n-1}{n+2}\a-1\right),\frac{2}{n-2}\right)$. By the continuity of \( V_c \) with respect to the parameter \( d \), there exists a $d=d(n,\a)\in(0,d_0)$ such that   $V_c>0$. It is obvious that for these fixed $\b,d$, we see that $U_c,V_c,W_c$ have positive lower bounds and the bounds only depend on $n,\a$. Therefore, there exist $L=L(n,\a)>0$, $M=M(n,\a)>0$,  $V_0=V_0(n,\a)>0$ and $W_0=W_0(n,\a)>0$ such that
		\begin{eqnarray}
			V&\ge&V_0-M\chi_{\{x\in B(x_0,2R): u(x)< L\e\}},\nonumber\\
			W&\ge&W_0-M\chi_{\{x\in B(x_0,2R): u(x)< L\e\}}.\nonumber
		\end{eqnarray}
		Then the proof is complete.

	\end{proof}

	\begin{lemma}\label{ac}
		Let $n\in[1,4)$, $\e>0$, $\g=1$, $\a\in[\frac{n+3}{n-1},p_S(n))$. Then the functions $U$, $V$, $W$ are defined as in Lemma \ref{k2}. Then there exist constants $\b=\b(n,\a)\in(0,\infty)$ if $n\in[1,2]$ or  $\b=\b(n,\a)\in(0,\frac{2}{n-2})$ if $n\in(2,4)$, $d=d(n,\a)>0$, $L=L(n,\a)>0$ and $M=M(n,\a)>0$ such that for any $\e>0$,
		\begin{equation}
			U\ge U_0>0,\nonumber
		\end{equation}
		\begin{equation}
			V\ge V_0-M\chi_{\{x\in B(x_0,2R): u(x)< L\e\}},\nonumber
		\end{equation}
		\begin{equation}
			W\ge W_0-M\chi_{\{x\in B(x_0,2R): u(x)< L\e\}},\nonumber
		\end{equation}
		where $U_0,V_0,W_0$ are positive constants which only depend on $n, \a$.	
		
	\end{lemma}
	
	\begin{proof}
		For $\g=1$, 	if $n\le 2$ and $\b>0$ or $n\in(2,4)$ and $\b\in(0,\frac{2}{n-2})$, we have
		\begin{equation}
			U=U_c=\left(\frac{2}{n}\left(1+\frac{1}{\b}\right)-1\right)\left(1+\frac{1}{\b}\right)>0,\nonumber
		\end{equation}
		Moreover, 
		\begin{eqnarray}
			V
			&=&V_c-\left(\b+2+\frac{4}{n}\left(1+\b\right)\right)\frac{\e}{u+\e}\nonumber\\
			&&+d\left(\frac{\left(\a-1\right)\left(\a-2\right)}{\b^2}\frac{2u\e+\e^2}{u^2}+\frac{2\left(\a-1\right)}{\b}\left(\frac{1}{\b}-1\right)\frac{\e}{u}\right),\nonumber\\
			W&=&W_c-\frac{d\b\e}{u+\e}-\frac{2\b^2}{n}\frac{\e^2+2u\e}{(u+\e)^2},\nonumber
		\end{eqnarray}
		where
		\begin{eqnarray}
			V_c&=&\frac{4}{n}(1+\b)+2\left(1-\a\right)+\b+\frac{d}{\b^2}\left(\a-\b-1\right)\left(\a-\b\right),\\
			W_c&=&\frac{2}{n}\b^2+d\left(\b+1-\a\right).
		\end{eqnarray}
	Before proceeding further with the discussion, we mention that the term in the second line of the expression for \( V \) $$\frac{2\left(\a-1\right)}{\b}\left(\frac{1}{\b}-1\right)\frac{\e}{u}$$ may be negative and may have no lower bound when \( \frac{\e}{u} \) is very large. So, the condition \( \a \ge \frac{n+3}{n-1}> 2 \) ($n<4$) is precisely to allow the term $\frac{\left(\a-1\right)\left(\a-2\right)}{\b^2}\frac{2u\e+\e^2}{u^2}$ to absorb this potential negative term in such cases.
		Now, we will divide the remaining argument into two cases for discussion.
		
		Case 1: $n\in[1,2]$. We set $\beta=\alpha-2\ge 3$, and take $d=d(\alpha)\in(0,1)$ sufficiently small. One readily sees that $U_c,V_c,W_c>0$, and these quantities depend only on $n,\alpha$. The condition $\alpha\ge\frac{n+3}{n-1}\ge 5$ ensures
		
		\begin{eqnarray}
			&&-\left(\b+2+\frac{4}{n}\left(1+\b\right)\right)\frac{\e}{u+\e}\nonumber\\
			&&+d\left(\frac{\left(\a-1\right)\left(\a-2\right)}{\b^2}\frac{2u\e+\e^2}{u^2}+\frac{2\left(\a-1\right)}{\b}\left(\frac{1}{\b}-1\right)\frac{\e}{u}\right)\nonumber\\
			&\ge& -\left(\b+2+\frac{4}{n}\left(1+\b\right)\right)\frac{\e}{u+\e}\nonumber\\
			&&+d\left(\frac{\left(\a-1\right)\left(\a-2\right)}{\b^2}\frac{2u\e+\e^2}{u^2}-\frac{2(\a-1)}{\b}\frac{\e}{u}\right)\nonumber\\
			&\ge&-10\a-\frac{d(\a-1)}{\a-2}\ge-10\a-\frac{4}{3}.
		\end{eqnarray}
		Simultaneously, we see that there exists a $L=L(n, \a)>0$ such that 
		\begin{eqnarray}
			&&-\left(\b+2+\frac{4}{n}\left(1+\b\right)\right)\frac{\e}{u+\e}\nonumber\\
			&&+d\left(\frac{\left(\a-1\right)\left(\a-2\right)}{\b^2}\frac{2u\e+\e^2}{u^2}-\frac{2}{\b}\frac{\e}{u}\right)\ge-\frac{V_c}{2},\nonumber\\
			&&-\frac{d\b\e}{u+\e}-\frac{2\b^2}{n}\frac{\e^2+2u\e}{(u+\e)^2}\ge-\frac{W_c}{2}.
		\end{eqnarray}
		on the set $\{x\in B(x_0,2R):u(x)\ge L\e\}$.
		
		Therefore, there exist $L=L(n,\a)>0$ and $M=M(\a)>0$ such that
		\begin{eqnarray}
			V&\ge&\frac{V_c}{2}-M\chi_{\{x\in B(x_0,2R): u(x)< L\e\}},\nonumber\\
			W&\ge&\frac{W_c}{2}-M\chi_{\{x\in B(x_0,2R): u(x)< L\e\}}.\nonumber
		\end{eqnarray}

			Case 2:	$n\in(2,4)$ and $\a\le\frac{n}{n-2}$. We first fix
			
			\[
			\beta  \in  \left( {{\left( \frac{1}{2} + \frac{2}{n}\right) }^{-1}\left( {\alpha  - 1 - \frac{2}{n}}\right) ,\alpha  - 1}\right) ,
			\]
			and then choose \(d = d\left( {n,\alpha }\right)  > 0\) (sufficiently small) such that \({U}_{c},{V}_{c},{W}_{c} > 0\) and they depend only on \(n,\alpha\). The remaining argument is similar to that in Case 1.

		Case 3: \(n \in  \left( {2,4}\right)\) and \(\alpha  \in  \left\lbrack  {\max \left( {\frac{n + 3}{n - 1},\frac{n}{n - 2}}\right) ,\frac{n + 2}{n - 2}}\right)\) . The same argument as in Lemma 3.5 can yield that there exist positive \(\beta ,d\) that depend only on \(n\) and \(\alpha\) such that \({U}_{c},{V}_{c},{W}_{c}\) have positive lower bounds and the bounds only depend on \(n,\alpha\) . The remaining argument is similar to that in Case 1 and so we omit these redundant details.
		
		Then the proof is complete.

	\end{proof}
	
	\begin{remark}\label{r37}
		In fact, the consideration of Case 2 in Lemma \ref{ac} is redundant because $\frac{n+3}{n-1}=\frac{n}{n-2}=3$ when $n=3$. However, in the proof of Lemma \ref{ab} and Lemma \ref{ac}, we have not required \( n \) to be a positive integer throughout, which lays the groundwork for our generalization to the general curvature-dimension condition in Section \ref{s5}.
	\end{remark}

	Before concluding this section, we present a useful remark. When \(\alpha\) is in some left neighborhood of \(\frac{n + 2}{n - 2}\) , we can provide more precise expressions for the coefficients \({U}_{0},{V}_{0},{W}_{0},L\left( {n,\alpha }\right)\) and \(M\left( {n,\alpha }\right)\) in Lemmas \ref{ab} and \ref{ac}.

	\begin{remark}\label{r31}
	
	For \(n \geq  3\) and \(\alpha  \in  \left\lbrack  {\frac{n + 3}{n - 1},\frac{n + 2}{n - 2}}\right)\) , we set
	
	\[
	\delta  = \frac{2}{n} - \frac{d\left( {\alpha  - 1 - \beta }\right) }{{\beta }^{2}}
	\]
	
	\[
	\theta  = \delta \left( {\alpha  - \beta }\right) .
	\]
		Let \({U}_{c},{V}_{c},{W}_{c}\) be defined by \eqref{uc},\eqref{vc} and \eqref{wc} respectively and now we fix
	
	\[
	\theta  = \frac{n - 1}{2n}\left( {\frac{n + 2}{n - 2} - \alpha }\right) 
	\]
	and
	
	\[
	\beta  = 2\left( {\frac{n - 1}{n + 2}\alpha  - 1}\right)  + \frac{2n\theta }{n + 2}
	\]
	
	\[
	= \frac{n - 1}{n + 2}\alpha  - \frac{n - 3}{n - 2} \in  \left( {\frac{4}{n + 2},\min \left( {\alpha  - 1,\frac{2}{n - 2}}\right) }\right) .
	\]
	Then direct computation yields that
	
	\[
	{U}_{c} = \frac{\left( {n - 1}\right) \left( {n - 2}\right) }{n\left( {n + 2}\right) }\left( {\frac{n + 2}{n - 2} - \alpha }\right) \frac{1 + \beta }{\beta^2}
	\]
	
	\[
	{V}_{c} = \theta  = \frac{n - 1}{2n}\left( {\frac{n + 2}{n - 2} - \alpha }\right) ,
	\]
	
	\[
	{W}_{c} = \delta {\beta }^{2} = \frac{{\beta }^{2}}{\alpha  - \beta } \times  \frac{n - 1}{2n}\left( {\frac{n + 2}{n - 2} - \alpha }\right) .
	\]
	From the expressions of \({U}_{c},{V}_{c}\) , and \({W}_{c}\) , we know that there exists a constant \(C = C\left( n\right)  > 0\) such that
	
	\[
	\min \left( {{U}_{c},{V}_{c},{W}_{c}}\right)  \geq  C\left( {\frac{n + 2}{n - 2} - \alpha }\right) .
	\]
	
	Define
	
	\[
	{U}_{0} = {U}_{c},\;{V}_{0} = \frac{{V}_{c}}{2},\;{W}_{0} = \frac{{W}_{c}}{2}.
	\]
	Next, we proceed with the analysis by dividing it into two cases.
	
	Case 1: \(n \geq  4\) . We set
	
	\[
	M = M\left( {n,\alpha }\right)  = \max \left( {1,d}\right) \beta ,\;L = L\left( {n,\alpha }\right)  = \frac{2\max \left( {1,d}\right) \beta }{\min \left( {{V}_{c},{W}_{c}}\right) }.
	\]
	Notice that
	
	\[
	d \leq  \frac{2{\beta }^{2}}{n\left( {\alpha  - 1 - \beta }\right) } = \frac{\frac{2}{n}{\beta }^{2}\left( {n - 2}\right) }{\frac{3\left( {n - 2}\right) }{n + 2}\alpha  - 1} \leq  \frac{6}{\frac{3\left( {n - 2}\right) }{n + 2} \times  \frac{n + 3}{n - 1} - 1} \leq  {12}.
	\]
	
	Thus, \(M\) has a uniform upper bound of 36, and \(L\) satisfies
	
	\[
	L \leq  C\left( n\right) {\left( \frac{n + 2}{n - 2} - \alpha \right) }^{-1}.
	\]
	
	Case 2: \(n \in  \lbrack 3,4)\) . From the proof of Lemma \ref{ac}, we set
	
	\[
	M = M\left( {n,\alpha }\right)  = \left( {\max \left( {1,d}\right)  + \frac{2}{n}\beta }\right) \beta  + 2 + \frac{4}{n}\left( {1 + \beta }\right)  + \frac{d}{\left( {\alpha  - 1}\right) \left( {\alpha  - 2}\right) }
	\]
	and
	
	\[
	L = \frac{2\left( {\beta  + 2 + \frac{4}{n}\left( {1 + \beta }\right) }\right) }{\min \left( {{V}_{c},{W}_{c}}\right) }.
	\]
	As in Case 1, we have
	
	\[
	d \leq  {12},\;\beta  \leq  3,\;M \leq  {100},\;L \leq  C\left( n\right) {\left( \frac{n + 2}{n - 2} - \alpha \right) }^{-1}.
	\]
	
	To summarize, when \(n \geq  3\) and \(\alpha  \in  \left\lbrack  {\frac{n + 3}{n - 1},\frac{n + 2}{n - 2}}\right)\) , by choosing appropriate \(d\) and \(\beta\) , the constants in Lemmas \ref{ab} and \ref{ac} satisfy the following estimates:
	
	\[
	\min \left( {{U}_{0},{V}_{0},{W}_{0}}\right)  \geq  C\left( n\right) \left( {\frac{n + 2}{n - 2} - \alpha }\right) ,\;M \leq  {100},\;L \leq  C\left( n\right) {\left( \frac{n + 2}{n - 2} - \alpha \right) }^{-1}.
	\]

	\end{remark}

	\section{Proof of Theorem \ref{mthm}}\label{s4}

	We are now ready to prove Theorem \ref{mthm}.  
	
	\begin{proof}[Proof of Theorem \ref{mthm}]
Without loss of generality, we assume the solution $u$ is defined on $B(x_0,4R)$, and we derive estimate \eqref{strong1} on $B(x_0,R)$. Otherwise, we apply this estimate to $B(y,\tfrac{R}{4})$ for $y\in B(x_0,R)$ and obtain the same estimate.

	According to Lemmas \ref{ps} and \ref{cl2}, to complete the proof of Theorem \ref{mthm}, we only need to prove the case where \(\alpha  \in  \left\lbrack  {\frac{n + 3}{n - 1},{p}_{S}\left( n\right) }\right)\) with \(n > 1\) . Specifically, we divide it into two cases.
	
	Case 1: \(n \geq  4\) and \(\alpha  \in  \left\lbrack  {\frac{n + 3}{n - 1},{p}_{S}\left( n\right) }\right)\) . Let \(u\) be a positive solution of \eqref{LE} on \(B\left( {{x}_{0},{2R}}\right)  \subset  {\mathcal{M}}^{n}\) . For an undetermined \(\varepsilon  > 0\) , we define the functions \(w\) by \eqref{tran} and \(F\) by \eqref{F} with \(\gamma  = 1\) and the \(\beta ,d\) for which Lemma \ref{ab} holds. By Lemma \ref{k1} and Lemma \ref{ab}, we have
	\begin{eqnarray}\label{41}
	{\Delta F} &\geq&  C{\left( u + \varepsilon \right) }^{\beta }{F}^{2} - {2KF} - {M}_{0}{u}^{\alpha  - 1}{\chi }_{\left\{  x \in  B\left( {x}_{0},2R\right)  : u\left( x\right)  < L\e \right\}  }F\nonumber\\
	&&+ 2\left( {\frac{1}{\beta } - 1 + \frac{u}{u + \varepsilon }}\right) \langle \nabla F,\nabla \ln w\rangle\nonumber\\
	&\geq&  C{\left( u + \varepsilon \right) }^{\beta }{F}^{2} - {2KF} - {M}_{0}{\left( L\varepsilon \right) }^{\alpha  - 1}F\nonumber\\
	&&+ 2\left( {\frac{1}{\beta } - 1 + \frac{u}{u + \varepsilon }}\right) \langle \nabla F,\nabla \ln w\rangle
	\end{eqnarray}
	on \(B\left( {{x}_{0},{2R}}\right)\) , where \(C = C\left( {n,\alpha }\right)  > 0,{M}_{0} = \max \left( {1,\frac{1}{d}}\right) M\left( {n,\alpha }\right) ,M=M\left( {n,\alpha }\right)\) and \(L = L\left( {n,\alpha }\right)\) are the same as those in Lemma \ref{ab}.
	
	Next, we derive the boundedness of \(F\) via \eqref{41} as follows. Let \(A = {\Phi F}\) be the auxiliary function, where \(\Phi\) is a cut-off function as in Lemma \ref{cut off}. Without loss of generality, we may assume that \(A\) attains its maximum inside \(B\left( {{x}_{0},\frac{3}{2}R}\right)\) or the proof is finished. We also assume that the maximal point \(x\) lies outside the cut locus of \({x}_{0}\) as before. By the chain rule, at point \(x\) , we see
	\begin{equation}\label{420}
			0 \geq  {\Delta A} = {F\Delta \Phi } + {\Phi \Delta F} - {2A}\frac{{\left| \nabla \Phi \right| }^{2}}{{\Phi }^{2}}. 
	\end{equation}
	By \eqref{41}, \eqref{420} and chain rule again, we have
	\begin{eqnarray}\label{43}
	0 &\geq&  {C\Phi }{\left( u + \varepsilon \right) }^{\beta }{F}^{2} - {2KA} - {M}_{0}{\left( L\varepsilon \right) }^{\alpha  - 1}A\nonumber\\
	&&- 2\left( {\frac{1}{\beta } - 1 + \frac{u}{u + \varepsilon }}\right) F\langle \nabla \Phi ,\nabla \ln w\rangle  + {F\Delta \Phi } - {2A}\frac{{\left| \nabla \Phi \right| }^{2}}{{\Phi }^{2}}. 
	\end{eqnarray}
	By the Cauchy-Schwarz inequality, we know that
	\begin{eqnarray}\label{44}
			&&- 2\left( {\frac{1}{\beta } - 1 + \frac{u}{u + \varepsilon }}\right) F\langle \nabla \Phi ,\nabla \ln w\rangle\nonumber\\
			&\geq&   - \frac{2}{C{\beta }^{2}}\frac{{\left| \nabla \Phi \right| }^{2}}{\Phi }F - \frac{C}{2}\Phi \frac{{\left| \nabla w\right| }^{2}}{{w}^{2}}F\nonumber\\
			&\geq&   - \frac{2}{C{\beta }^{2}}\frac{{\left| \nabla \Phi \right| }^{2}}{\Phi }F - \frac{C}{2}\Phi {\left( u + \varepsilon \right) }^{\beta }{F}^{2} 
	\end{eqnarray}

	Substituting \eqref{44} into \eqref{43} yields
\begin{equation}\label{45}
		0 \geq  \frac{C}{2}\Phi {\left( u + \varepsilon \right) }^{\beta }{F}^{2} - {2KA} - {M}_{0}{\left( L\varepsilon \right) }^{\alpha  - 1}A + {F\Delta \Phi } - 2\left( {\frac{1}{C{\beta }^{2}} + 1}\right) \frac{{\left| \nabla \Phi \right| }^{2}}{{\Phi }^{2}}A.
\end{equation}
	Multiplying \(\Phi \left( x\right)\) on both sides of \eqref{45} and using the property of \(\Phi\) , we have
	
	\[
	A\left( x\right)  \leq  C\left( {n,\alpha }\right) {\left( u + \varepsilon \right) }^{-\beta }\left( x\right) \left( {K + \frac{1}{{R}^{2}} + {\left( L\varepsilon \right) }^{\alpha  - 1}}\right)
	\]
	
	\[
	\leq  C\left( {n,\alpha }\right) {\varepsilon }^{-\beta }\left( {K + \frac{1}{{R}^{2}} + {\left( L\varepsilon \right) }^{\alpha  - 1}}\right) ,
	\]
	
	which implies
	\begin{equation}\label{46}
		\mathop{\sup }\limits_{{B\left( {{x}_{0},R}\right) }}{\left( u + \varepsilon \right) }^{-\beta }\left( {\frac{{\left| \nabla u\right| }^{2}}{{u}^{2}} + {u}^{\alpha  - 1}}\right)  \leq  C\left( {n,\alpha }\right) {\varepsilon }^{-\beta }\left( {K + \frac{1}{{R}^{2}} + {\left( L\varepsilon \right) }^{\alpha  - 1}}\right).
	\end{equation}
	
	Now, we choose \(\varepsilon  > 0\) such that
	\begin{equation}\label{47}
		{\left( L\varepsilon \right) }^{\alpha  - 1} = K + \frac{1}{{R}^{2}}.
	\end{equation}
	For \(y \in  B\left( {{x}_{0},R}\right)\), if \(u\left( y\right)  \leq  \varepsilon\), by \eqref{46} and \eqref{47}, we see
	\begin{equation}\label{48}
			\left( {\frac{{\left| \nabla u\right| }^{2}}{{u}^{2}} + {u}^{\alpha  - 1}}\right) \left( y\right)  \leq  C\left( {n,\alpha }\right) \left( {K + \frac{1}{{R}^{2}}}\right).
	\end{equation}
	If \(u\left( y\right)  \geq  \varepsilon\), \eqref{46} yields that \(u\left( y\right)  \leq  C\left( {n,\alpha }\right) \varepsilon\) because 
	\[
	\alpha  - 1 - \beta  > \frac{n + 3}{n - 1} - 1 - \frac{2}{n - 2} > 0.
	\] Further, we also obtain \eqref{48} by \eqref{46} and \eqref{47}. Therefore,
	
	\[
	\mathop{\sup }\limits_{{B\left( {{x}_{0},R}\right) }}\left( {\frac{{\left| \nabla u\right| }^{2}}{{u}^{2}} + {u}^{\alpha  - 1}}\right)  \leq  C\left( {n,\alpha }\right) \left( {K + \frac{1}{{R}^{2}}}\right) .
	\]
	Thus, we have completed the proof of Theorem \ref{mthm} for \(n \geq  4\) .

	Case 2: \(n \in  \left( {1,4}\right)\) and \(\alpha  \in  \left\lbrack  {\frac{n + 3}{n - 1},{p}_{S}\left( n\right) }\right)\). Let \(u\) be a positive solution of (1.2) on \(B\left( {{x}_{0},{2R}}\right)\). For an undetermined \(\varepsilon  > 0\), we define the functions \(w\) by \eqref{tran2} and \(G\) by \eqref{G} with \(\gamma  = 1\) and the \(\beta ,d\) for which Lemma \ref{ac} holds. By Lemma \ref{k2} and Lemma \ref{ac}, we have
	\begin{equation}\label{49}
		{\Delta G} \geq  C{w}^{-1}{G}^{2} - {2KG} + \frac{2}{\beta }\langle \nabla G,\nabla \ln w\rangle  - {M}_{0}{\left( L\varepsilon \right) }^{\alpha  - 1}G
	\end{equation}
	on \(B\left( {{x}_{0},{2R}}\right)\), where \(C = C\left( {n,\alpha }\right)  > 0,{M}_{0} = \max \left( {1,\frac{1}{d}}\right) M\left( {n,\alpha }\right) ,M=M\left( {n,\alpha }\right)\) and \(L = L\left( {n,\alpha }\right)\) are the same as those in Lemma \ref{ac}.
	
	Next, we derive the boundedness of \(G\) via \eqref{49} as follows. Let \(A = {\Phi G}\) be the auxiliary function, where \(\Phi\) is a cut-off function as in Lemma \ref{cut off}. Without loss of generality, we may assume that \(A\) attains its maximum inside \(B\left( {{x}_{0},\frac{3}{2}R}\right)\) or the proof is finished. We also assume that the maximal point \(x\) lies outside the cut locus of \({x}_{0}\) as before. By the chain rule, at point \(x\), we see
	\begin{equation}\label{410}
			0 \geq  {\Delta A} = {G\Delta \Phi } + {\Phi \Delta G} - {2A}\frac{{\left| \nabla \Phi \right| }^{2}}{{\Phi }^{2}}.
	\end{equation}
	
	By \eqref{49}, \eqref{410} and chain rule, we have
	\begin{eqnarray}\label{411}
			0 &\geq&  {C\Phi }{w}^{-1}{G}^{2} - {2KA} - {M}_{0}{\left( L\varepsilon \right) }^{\alpha  - 1}A\nonumber\\
			&&- \frac{2}{\beta }G\langle \nabla \Phi ,\nabla \ln w\rangle  + {G\Delta \Phi } - {2A}\frac{{\left| \nabla \Phi \right| }^{2}}{{\Phi }^{2}}.
	\end{eqnarray}
	The Cauchy-Schwarz inequality yields
	\begin{eqnarray}\label{412}
		- \frac{2}{\beta }G\langle \nabla \Phi ,\nabla \ln w\rangle  &\geq&   - \frac{2}{C{\beta }^{2}}\frac{{\left| \nabla \Phi \right| }^{2}}{\Phi }G - \frac{C}{2}\Phi \frac{{\left| \nabla w\right| }^{2}}{{w}^{2}}G\nonumber\\
		&\geq&   - \frac{2}{C{\beta }^{2}}\frac{{\left| \nabla \Phi \right| }^{2}}{\Phi }G - \frac{C}{2}\Phi {\left( u + \varepsilon \right) }^{\beta }{G}^{2}.
	\end{eqnarray}
	Substituting \eqref{412} into \eqref{411} first, then using the property of \(\Phi\), we have
	
	\[
	A\left( x\right)  \leq  C\left( {n,\alpha }\right) {\left( u + \varepsilon \right) }^{-\beta }\left( x\right) \left( {K + \frac{1}{{R}^{2}} + {\left( L\varepsilon \right) }^{\alpha  - 1}}\right)
	\]
	
	\[
	\leq  C\left( {n,\alpha }\right) {\varepsilon }^{-\beta }\left( {K + \frac{1}{{R}^{2}} + {\left( L\varepsilon \right) }^{\alpha  - 1}}\right) ,
	\]
	
	which implies
	\begin{equation}\label{413}
			\mathop{\sup }\limits_{{B\left( {{x}_{0},R}\right) }}{\left( u + \varepsilon \right) }^{-\beta }\left( {\frac{{\left| \nabla u\right| }^{2}}{{\left( u + \varepsilon \right) }^{2}} + {u}^{\alpha  - 1}}\right)  \leq  C\left( {n,\alpha }\right) {\varepsilon }^{-\beta }\left( {K + \frac{1}{{R}^{2}} + {\left( L\varepsilon \right) }^{\alpha  - 1}}\right) . 
	\end{equation}

	Now, we choose \(\varepsilon  > 0\) such that \eqref{47} holds. For \(y \in  B\left( {{x}_{0},R}\right)\) , if \(u\left( y\right)  \leq  \varepsilon\) , by \eqref{413} and \eqref{47}, we see
	\begin{equation}\label{414}
		{u}^{\alpha  - 1}\left( y\right)  \leq  C\left( {n,\alpha }\right) \left( {K + \frac{1}{{R}^{2}}}\right).
	\end{equation}
	If \(u\left( y\right)  \geq  \varepsilon\), \eqref{46} yields that \(u\left( y\right)  \leq  C\left( {n,\alpha }\right) \varepsilon\) because \(\alpha  - 1 - \beta  > 0\). Further, we also obtain \eqref{414} by \eqref{413} and \eqref{47}. Therefore,
\begin{equation}\label{415}
	\mathop{\sup }\limits_{{B\left( {{x}_{0},R}\right) }}{u}^{\alpha  - 1} \leq  C\left( {n,\alpha }\right) \left( {K + \frac{1}{{R}^{2}}}\right) .
\end{equation}
		Combining \eqref{415} and Lemma \ref{cl1}, we have the desired estimate in Theorem \ref{mthm} and thus complete the proof of Theorem \ref{mthm}.

	\end{proof} 
   	
	\begin{remark}
		
		For \(n \geq  3\) and \(\alpha  \in  \left\lbrack  {\frac{n + 3}{n - 1},\frac{n + 2}{n - 2}}\right)\), we will precisely compute the constant \(C\left( {n,\alpha }\right)\) in the estimate \eqref{strong1}. From Remark \ref{r31}, we know that
		
		\[
		0 < c\left( n\right)  \leq  d \leq  {12},\;\beta  \in  \left( {\frac{4}{n + 2},\min \left( {\alpha  - 1,\frac{2}{n - 2}}\right) }\right)
		\]
		and
		
		\[
		\min \left( {{U}_{0},{V}_{0},{W}_{0}}\right)  \geq  C\left( n\right) \left( {\frac{n + 2}{n - 2} - \alpha }\right) ,\;M \leq  {100},\;L \leq  C\left( n\right) {\left( \frac{n + 2}{n - 2} - \alpha \right) }^{-1}.
		\]
		Substituting these more refined estimates into the proof of Theorem \ref{mthm}, we can easily obtain that the constant \(C\left( {n,\alpha }\right)\) in \eqref{46} and \eqref{413} can be expressed as
		
		\[
		C\left( n\right) {\left( \frac{n + 2}{n - 2} - \alpha \right) }^{-2}.
		\]
		Furthermore, the constant in the final estimate can be expressed as
		
		\[
		C\left( n\right) {\left( \frac{n + 2}{n - 2} - \alpha \right) }^{-\left( {2 + \frac{\left( {\alpha  + 1}\right) \beta }{\alpha  - 1 - \beta }}\right) }.
		\]
		It should be noted that at this point, \(2 + \frac{\left( {\alpha  + 1}\right) \beta }{\alpha  - 1 - \beta } \leq  \theta \left( n\right)\) , where \(\theta \left( n\right)\) is a positive constant depending only on \(n\). So, we have proven the claim in Remark \ref{r11}.

	\end{remark}

	\section{Generalization under Bakry-\'{E}mery curvature}\label{s5}

\subsection{A generalization of Theorem \ref{mthm}}
In this subsection, we generalize previous results for equation \eqref{LE} to the following equation with gradient term:
\begin{equation}\label{vle}
	{\Delta }_{V}u + {u}^{\alpha } = 0,
\end{equation}where \({\Delta }_{V}u \mathrel{\text{ := }} {\Delta u} + \langle \nabla u,V\rangle\) and \(V\) is a \({C}^{1}\) vector field. When \(V\) vanishes, then \eqref{vle} becomes \eqref{LE}.

To derive universal bounds and gradient estimates, we need to recall some previous lemmas and their proofs and analyse how to generalize them for equation \eqref{vle}. First, we need analogues of Lemma \ref{k1} and Lemma \ref{k2} which yield elliptic inequalities of the first and second type auxiliary functions of solutions respectively. For a brief description, we define \(\infty\) -dimensional and \(m\) -dimensional Bakry-Émery Ricci curvature respectively.

\[
{\operatorname{Ric}}_{V} \mathrel{\text{ := }} \operatorname{Ric} - \frac{1}{2}{\mathcal{L}}_{V}g
\]

\[
{\operatorname{Ric}}_{V}^{m} \mathrel{\text{ := }} {\operatorname{Ric}}_{V} - \frac{1}{m - n}{V}^{\flat } \otimes  {V}^{\flat },
\]
where \(m > n\) and \(n\) is the dimension of the manifold, \({\mathcal{L}}_{V}g\) is the Lie derivative of Riemannian metric \(g\) on the vector field \(V\) and \({V}^{\flat }\) is the dual 1-form of \(V\) under metric \(g\). From the classical Bochner formula, one can derive its generalization:
\begin{equation}\label{vbo}
	{\Delta }_{V}{\left| \nabla w\right| }^{2} = 2{\left| {\nabla }^{2}w\right| }^{2} + 2\left\langle  {\nabla {\Delta }_{V}w,\nabla w}\right\rangle   + 2{\operatorname{Ric}}_{V}\left( {\nabla w,\nabla w}\right) , 
\end{equation}where \(w \in  {C}^{3}\left( \mathcal{M}^n\right)\). We have the following Leibniz rule.\begin{equation}\label{vleib}
{\Delta }_{V}\left( {uv}\right)  = v{\Delta }_{V}u + u{\Delta }_{V}v + 2\langle \nabla u,\nabla v\rangle,
\end{equation}
where \(u,v \in  {C}^{2}\left( {\mathcal{M}}^{n}\right)\) . We also need the following inequality to substitute \eqref{id}:
\begin{eqnarray}\label{54}
	{\left| {\nabla }^{2}w\right| }^{2} &\geq&  \frac{1}{n}{\left( \Delta w\right) }^{2} = \frac{1}{n}\left( {{\Delta }_{V}w-\langle \nabla w,V\rangle }\right)\nonumber\\
	&\geq&  \frac{1}{n}\left( {\frac{{\left( {\Delta }_{V}w\right) }^{2}}{\frac{m}{n}} - \frac{{\left( \langle \nabla w,V\rangle \right) }^{2}}{\frac{m}{n} - 1}}\right)\nonumber\\
	&=& \frac{1}{m}{\left( {\Delta }_{V}w\right) }^{2} - \frac{1}{m - n}{V}^{\flat } \otimes  {V}^{\flat }\left( {\nabla w,\nabla w}\right),
\end{eqnarray}
where \(m > n\) is arbitrary and \(w \in  {C}^{2}\left( \mathcal{M}\right)\). Substituting \eqref{54} into \eqref{vbo} yields the generalized Bochner inequality:
\begin{equation}\label{55}
	{\Delta }_{V}{\left| \nabla w\right| }^{2} \geq  \frac{2}{m}{\left( {\Delta }_{V}w\right) }^{2} + 2\left\langle  {\nabla {\Delta }_{V}w,\nabla w}\right\rangle   + 2{\operatorname{Ric}}_{V}^{m}\left( {\nabla w,\nabla w}\right).
\end{equation}

Now, for a positive solution \(u\) in \(B\left( {{x}_{0},{2R}}\right)\) , we also define \(w\) and \(F\) by \eqref{tran} and \eqref{F} and call \(F\) the first type auxiliary function, and define \(w\) and \(G\) by \eqref{tran2} and \eqref{G} and call \(G\) the second type auxiliary function.

By substituting Leibniz rule and Bochner formula with \eqref{vleib} and \eqref{55} in the proof of Lemma \ref{k1} (inequality \eqref{55} substitutes the original Bochner formula and \eqref{id}), through exactly the same calculation process as in Lemma \ref{k1}, we have the following generalization of Lemma \ref{k1}.

\begin{lem}\label{vk1}
 Let \(u\) be a positive solution of \eqref{vle} on \(B\left( {{x}_{0},{2R}}\right)  \subset  {\mathcal{M}}^{n}\) , \(w\) be defined by \eqref{tran} with \(\beta  \neq  0\) and \(F\) be its first type auxiliary function as \eqref{F} with \(\gamma ,\varepsilon ,d \geq  0\) . Then for any \(m > n\) , we have

\begin{eqnarray}
	{\left( u + \varepsilon \right) }^{\beta \gamma }{\Delta }_{V}F &\geq&  2{w}^{-2}{\operatorname{Ric}}_{V}^{m}\left( {\nabla w,\nabla w}\right)\nonumber\\
	&&+ 2\left( {\frac{1}{\beta } - 1 + \gamma \frac{u}{u + \varepsilon }}\right) {\left( u + \varepsilon \right) }^{\beta \gamma }\langle \nabla F,\nabla \ln w\rangle\nonumber\\
	&&+ U\frac{{\left| \nabla w\right| }^{4}}{{w}^{4}} + V\frac{{\left| \nabla w\right| }^{2}}{{w}^{2}}{u}^{\alpha  - 1} + W{u}^{2\left( {\alpha  - 1}\right) }
\end{eqnarray}
on \(B\left( {{x}_{0},{2R}}\right)\), where

\begin{align*}
	U &= \frac{2}{m}{\left( 1 + \frac{1}{\beta }\right) }^{2} + \left( {\frac{\gamma }{\beta } - {\gamma }^{2}}\right) \frac{{u}^{2}}{{\left( u + \varepsilon \right) }^{2}} + 2\left( {1 - \frac{1}{\beta }}\right) \frac{\gamma u}{u + \varepsilon } - 2, \\
	V &= \frac{4}{m}\left( {1 + \beta }\right)  + 2\left( {1 - \alpha }\right)  + d\left( {\frac{\alpha \left( {\alpha  - 1}\right) }{{\beta }^{2}} - \frac{2\left( {\alpha  - 1}\right) }{\beta }}\right) \\
	&\quad + \frac{\gamma u}{u + \varepsilon }\left( {\beta  + d\left( {\left( {\frac{1}{\beta } - \gamma }\right) \frac{u}{u + \varepsilon } + 2 - \frac{2}{\beta }}\right) }\right) , \\
	W &= \frac{2{\beta }^{2}}{m} + d\left( {\frac{\beta \gamma u}{u + \varepsilon } + 1 - \alpha }\right) .
\end{align*}
\end{lem}

Similarly, the generalization of Lemma \ref{k2} is as follows:

\begin{lem}\label{vk2}
 Let \(u\) be a positive solution of \eqref{vle} on \(B\left( {{x}_{0},{2R}}\right)  \subset  {\mathcal{M}}^{n}\) , \(w\) be defined as \eqref{tran2} with \(\beta  \neq  0\) and \(G\) be its second type auxiliary function as \eqref{G} with \(\varepsilon ,d > 0\) and \(\gamma  = 1\) . Then for any \(m > n\) , we have
 \begin{eqnarray}
 	{w}^{-1}{\Delta }_{V}G \geq  2{w}^{-2}{\operatorname{Ric}}_{V}^{m}\left( {\nabla w,\nabla w}\right)  + \frac{2}{\beta }{w}^{-\gamma }\langle \nabla G,\nabla \ln w\rangle\nonumber\\
 	+ U\frac{{\left| \nabla w\right| }^{4}}{{w}^{4}} + V\frac{{\left| \nabla w\right| }^{2}}{{w}^{2}}{u}^{\alpha  - 1} + W{u}^{2\left( {\alpha  - 1}\right) }
 \end{eqnarray}
on \(B\left( {{x}_{0},{2R}}\right)\) , where
\begin{align*}
	U &= \left( {\frac{2}{m}\left( {1 + \frac{1}{\beta }}\right)  - 1}\right) \left( {1 + \frac{1}{\beta }}\right) , \\
	V &= \left( {\frac{4}{m}\left( {1 + \beta }\right)  + 2 + \beta }\right) \frac{u}{u + \varepsilon } - {2\alpha } + \frac{2d}{\beta }\left( {\frac{\alpha  - 1}{\beta }\frac{u + \varepsilon }{u} - 1}\right) \\
	&\quad + d\left( {\frac{\left( {\alpha  - 1}\right) \left( {\alpha  - 2}\right) }{{\beta }^{2}}\frac{{\left( u + \varepsilon \right) }^{2}}{{u}^{2}} + 1 + \frac{1}{\beta } - \frac{2\left( {\alpha  - 1}\right) }{\beta }\frac{u + \varepsilon }{u}}\right) , \\
	W &= \frac{2{\beta }^{2}}{m}{\left( \frac{u}{u + \varepsilon }\right) }^{2} + d\left( {\frac{\beta \gamma u}{u + \varepsilon } + 1 - \alpha }\right) .
\end{align*}
\end{lem}

Via the generalized Laplacian comparison theorem from Bakry-Qian \cite[Theorem 4.2]{BQ}, we also have a good cut-off function as in Lemma \ref{cut off}.

\begin{lem}\label{vcut}
	Let \(\left( {{\mathcal{M}}^{n},g}\right)\) be an \(n\)-dimensional complete Riemannian manifold with \({Ri}{c}_{V}^{m} \geq \; - {Kg}\) , where \(K \geq  0\) and \(m > n\). Then for any \(R > 0\) , we have a Lipschitz function \(\Phi\) on \(B\left( {{x}_{0},{2R}}\right)\) such that
	
	(i) \(\Phi \left( x\right)  = \phi \left( {d\left( {{x}_{0},x}\right) }\right)\), where \(d\left( {{x}_{0}, \cdot  }\right)\) is the distance function from \({x}_{0}\) and \(\phi\) is a non-increasing function on \(\lbrack 0,\infty )\) and
	
	\[
	\Phi \left( x\right)  = \left\{  \begin{array}{lll} 1 &\qquad \text{ if } & x \in  B\left( {{x}_{0},R}\right) \\  0 & \qquad\text{ if } & x \in  B\left( {{x}_{0},{2R}}\right)  \setminus  B\left( {{x}_{0},\frac{3}{2}R}\right) . \end{array}\right.
	\]
	
	(ii) On $\left\{  {x \in  B\left( {{x}_{0},{2R}}\right)  : \Phi \left( x\right)  > 0}\right\}$,
	
	\[
	\frac{\left| \nabla \Phi \right| }{{\Phi }^{\frac{1}{2}}} \leq  \frac{C}{R}
	\]
	
	(iii)
	
	\[
	{\Delta }_{V}\Phi  \geq   - \frac{C\sqrt{mK}}{R}\coth \left( {\sqrt{\frac{K}{m}}R}\right)  - \frac{C}{{R}^{2}} \geq  -C\left( m\right) \frac{\sqrt{K}R + 1}{{R}^{2}}
	\]
	holds on \(B\left( {{x}_{0},{2R}}\right)\) in the distribution sense and pointwise outside the cut locus of \({x}_{0}\). Here, \(C\) is a universal constant.
\end{lem}

Except for replacing the cut-off function in Lemma \ref{cut off} with that in Lemma \ref{vcut}, we can generalize all the lemmas in Section 3 to the current setting through exactly the same reasoning. We list the generalized results directly below without rewriting the proofs.

\begin{lem}
	Let \(\left( {{\mathcal{M}}^{n},g}\right)\) be an \(n\)-dimensional complete Riemannian manifold with \({Ri}{c}_{V}^{m} \geq \; - {Kg}\), where \(K \geq  0\) and \(m > n\) . If \(u\) is a positive solution of \eqref{vle} on \(B\left( {{x}_{0},{2R}}\right)  \subset \; {\mathcal{M}}^{n}\) with \(\alpha  < \frac{m + 3}{m - 1}\), then
	\begin{equation}
		\mathop{\sup }\limits_{{B\left( {{x}_{0},R}\right) }}\frac{{\left| \nabla u\right| }^{2}}{{u}^{2}} \leq  C\left( {m,\a}\right) \left( {K + \frac{1}{{R}^{2}}}\right).
	\end{equation}
\end{lem}

\begin{lem}
	Let \(\left( {{\mathcal{M}}^{n},g}\right)\) be an \(n\)-dimensional complete Riemannian manifold with \({Ri}{c}_{V}^{m} \geq \; - {Kg}\), where \(K \geq  0\) and \(m > n\). If \(u\) is a positive solution of \eqref{vle} on \(B\left( {{x}_{0},{2R}}\right)  \subset \; {\mathcal{M}}^{n}\) with \(\alpha  \in  \mathbb{R}\), then
\begin{equation}
	\mathop{\sup }\limits_{{B\left( {{x}_{0},R}\right) }}\frac{{\left| \nabla u\right| }^{2}}{{u}^{2}} \leq  C\left( {m,\alpha }\right) \left( {K + \frac{1}{{R}^{2}} + \mathop{\sup }\limits_{{B\left( {{x}_{0},{2R}}\right) }}{u}^{\alpha  - 1}}\right).
\end{equation}
	Furthermore, \(C\left( {m,\alpha }\right)  = C\left( m\right) \alpha\) when \(\alpha  > 1\).
\end{lem}

\begin{lem}
	Let \(\left( {{\mathcal{M}}^{n},g}\right)\) be an \(n\)-dimensional complete Riemannian manifold with \({\operatorname{Ric}}_{V}^{m} \geq \; - {Kg}\), where \(K \geq  0\) and \(m > n\). If \(u\) is a positive solution of \eqref{vle} on \(B\left( {{x}_{0},{2R}}\right)  \subset \; {\mathcal{M}}^{n}\) with \(\alpha  \in  \mathbb{R}\), then
\begin{equation}
	\mathop{\sup }\limits_{{B\left( {{x}_{0},R}\right) }}{u}^{\alpha  - 1} \leq  C\left( {m,\alpha }\right) \left( {K + \frac{1}{{R}^{2}} + \mathop{\sup }\limits_{{B\left( {{x}_{0},{2R}}\right) }}\frac{{\left| \nabla u\right| }^{2}}{{u}^{2}}}\right).
\end{equation}
	Furthermore, \(C\left( {m,\alpha }\right)  = C\left( m\right) {\alpha }^{2}\) when \(\alpha  > 1\).
\end{lem}

As pointed out in Remark \ref{r37}, in the proof of Lemma \ref{ab} and Lemma \ref{ac}, we have not required \(n\) to be a positive integer throughout. From Lemma \ref{vk1} and \ref{vk2}, we can easily see that the only difference between \(U,V,W\) and those in Lemmas \ref{k1} and \ref{k2} is that \(n\) is replaced by \(m\) . Thus, through exactly the same process as in Lemma \ref{ab} and \ref{ac}, we obtain the lower bound estimate of the leading coefficients.

\begin{lem}\label{vab}
	Let \(m \geq  4,\varepsilon  > 0,\gamma  = 1\) and \(\alpha  \in  \left\lbrack  {\frac{m + 3}{m - 1},{p}_{S}\left( m\right) }\right)\). The functions \(U,V,W\) are defined as in Lemma \ref{vk1}. Then there exist constants \(\beta  = \beta \left( {m,\alpha }\right)  \in \; \left( {0,\frac{2}{m - 2}}\right), d = d\left( {m,\alpha }\right) > 0, L = L\left( {m,\alpha }\right)  > 0\) and \(M = M\left( {m,\alpha }\right)  > 0\) such that for any \(\varepsilon  > 0\),
	
	\[
	U \geq  {U}_{0} > 0,
	\]
	
	\[
	V \geq  {V}_{0} - M{\chi }_{\left\{  x \in  B\left( {x}_{0},2R\right)  : u\left( x\right)  < L\e \right\}  },
	\]
	
	\[
	W \geq  {W}_{0} - M{\chi }_{\{ x \in  B\left( {{x}_{0},{2R}}\right)  : u\left( x\right)  < {L\e }\} },
	\]
	where \({U}_{0},{V}_{0},{W}_{0}\) are positive constants which only depend on \(m,\alpha\).
\end{lem}
\begin{lem}\label{vac}
	Let \(m \in  \left( {1,4}\right), \varepsilon  > 0, \gamma  = 1,\alpha  \in  \left\lbrack  {\frac{m + 3}{m - 1},{p}_{S}\left( m\right) }\right)\). Then the functions \(U\), \(V\), \(W\) are defined as in Lemma \ref{vk2}. Then there exist constants \(\beta  = \beta \left( {m,\alpha }\right)  \in \; \left( {0,\alpha  - 1}\right)\) if \(m \in  (1,2\rbrack\) or \(\beta  = \beta \left( {m,\alpha }\right)  \in  \left( {0,\min \left( {\frac{2}{m - 2},\alpha  - 1}\right) }\right)\) if \(m \in  \left( {2,4}\right)\) , \(d = d\left( {m,\alpha }\right)  > 0,L = L\left( {m,\alpha }\right)  > 0\) and \(M = M\left( {m,\alpha }\right)  > 0\) such that for any \(\varepsilon  > 0,\)
	
	\[
	U \geq  {U}_{0} > 0,
	\]
	
	\[
	V \geq  {V}_{0} - M{\chi }_{\left\{  x \in  B\left( {x}_{0},2R\right)  : u\left( x\right)  < L\e \right\}  },
	\]
	
	\[
	W \geq  {W}_{0} - M{\chi }_{\{ x \in  B\left( {{x}_{0},{2R}}\right)  : u\left( x\right)  < {L\e }\} },
	\]
		where \({U}_{0},{V}_{0},{W}_{0}\) are positive constants which only depend on \(m,\alpha\).
\end{lem}

Using Lemmas \ref{vab} and \ref{vac} and the same reasoning as in Section 4, we derive the following generalization of Theorem \ref{mthm}.

\begin{theorem}\label{vmthm}
	 Let \(\left( {{\mathcal{M}}^{n},g}\right)\) be an \(n\)-dimensional complete Riemannian manifold with \({\operatorname{Ric}}_{V}^{m} \geq \; - {Kg}\), where \(K \geq  0\) and \(m > n\). If \(u\) is a positive solution of \eqref{vle} on \(B\left( {{x}_{0},{2R}}\right)  \subset \; {\mathcal{M}}^{n}\) with \(\alpha  \in  \left( {-\infty ,{p}_{S}\left( m\right) }\right)\), then
		\[
	\mathop{\sup }\limits_{{B\left( {{x}_{0},R}\right) }}\left( {\frac{{\left| \nabla u\right| }^{2}}{{u}^{2}} + {u}^{\alpha  - 1}}\right)  \leq  C\left( {m,\alpha }\right) \left( {K + \frac{1}{{R}^{2}}}\right).
	\]
\end{theorem}

A natural corollary is the following Liouville theorem:

\begin{corollary}
	Let \(\left( {{\mathcal{M}}^{n},g}\right)\) be an $n$-dimensional complete Riemannian manifold with \({\operatorname{Ric}}_{V}^{m} \geq  0\). Then there are no positive solutions to equation \eqref{vle} on the manifold \({\mathcal{M}}^{n}\) with \(\alpha  \in  \left( {-\infty ,{p}_{S}\left( m\right) }\right)\). Furthermore, any non-negative solution to equation \eqref{vle} on \({\mathcal{M}}^{n}\) with \(\alpha  \in  \left( {0,{p}_{S}\left( m\right) }\right)\) is trivial.
\end{corollary}

For the case of Euclidean space, we can obtain the following estimate from Theorem \ref{vmthm}:

 \begin{theorem}
 	Let \(u\) be a positive solution of \eqref{vle} in \(\Omega  \subset  {\mathbb{R}}^{n}\) with \(\alpha  \in \; \left( {-\infty ,{p}_{S}\left( n\right) }\right)\). For \(x \in  \Omega\), we have
 		\[
 	\left( {\frac{{\left| \nabla u\right| }^{2}}{{u}^{2}} + {u}^{\alpha  - 1}}\right) \left( x\right)  \leq  C\left( {n,\alpha }\right) \left( {\frac{1}{{d}^{2}\left( {x,\partial \Omega }\right) } + \mathop{\max }\limits_{{B\left( {x,\frac{1}{2}d\left( {x,\partial \Omega }\right) }\right) }}\left( {{\left| V\right| }^{2} + \left| {\nabla V}\right| }\right) }\right)
 	\]
 	where \(d\left( {x,\partial \Omega }\right)\) is the distance from \(x\) to \(\partial \Omega\).
 \end{theorem}

\begin{proof}
	For \(x \in  \Omega\) , we define \(R = \frac{1}{8}d\left( {x,\partial \Omega }\right)\) and choose an undetermined cut-off function \({\varphi }_{x} \in  {C}_{0}^{\infty }\left( {B\left( {x,{4R}}\right) }\right)\) satisfying
	\begin{itemize}
		\item \({\varphi }_{x} = 1\) in \(B\left( {x,{2R}}\right)\) and \(0 \leq  {\varphi }_{x} \leq  1\),
		\item \(\left| {\nabla {\varphi }_{x}}\right|  \leq  C{R}^{-1}\) , where \(C = C\left( n\right)\).
	\end{itemize}
	Hence \(u\) solves equation \eqref{vle} on \(B\left( {x,{2R}}\right)\) with vector field \({\varphi }_{x}V\). Then we estimate the lower bound of \({\operatorname{Ric}}_{{\varphi }_{x}V}^{m}\), where \(m \in  \left( {n,\frac{2\left( {\alpha  + 1}\right) }{\alpha  - 1}}\right)\) if \(\alpha  > 1\) or \(m = {2n}\) if \(\alpha  \leq  1\).
	\begin{eqnarray}
		\operatorname{Ric}_{{\varphi }_{x}V}^{m} &=&  - \frac{1}{2}{{\mathcal{L}}_{{\varphi }_{x}V}g} - \frac{{\varphi }_{x}^{2}}{m - n}{V}^{\flat } \otimes  {V}^{\flat } \nonumber \\
		&=&  - \operatorname{Sym}\left( {{V}^{\flat } \otimes  d{\varphi }_{x}}\right)  - \frac{{\varphi }_{x}}{2}{{\mathcal{L}}_{V}g} - \frac{{\varphi }_{x}^{2}}{m - n}{V}^{\flat } \otimes  {V}^{\flat } \nonumber \\
		&\geq&   - C\left| {d{\varphi }_{x}}\right| \left| V\right| g - {\varphi }_{x}{\operatorname{Ric}}_{V}^{m} \nonumber \\
		&\geq&   - C{R}^{-1}\left| V\right| g - C{\varphi }_{x}\left( {\left| {\nabla V}\right|  + {\left| V\right| }^{2}}\right) g \nonumber \\
		&\geq&   - C\left( {{R}^{-2} + \mathop{\max }\limits_{{B\left( {x,{4R}}\right) }}\left( {{\left| V\right| }^{2} + \left| {\nabla V}\right| }\right) }\right) g, 
	\end{eqnarray}
		where \(\operatorname{Sym}\left( {{V}^{\flat } \otimes  d{\varphi }_{x}}\right)\) is the symmetrization of \({V}^{\flat } \otimes  d{\varphi }_{x}\) and \(g\) is the standard metric in \({\mathbb{R}}^{n}\). Then we finish the proof by applying Theorem \ref{vmthm}.
\end{proof}

When \(\Omega\) becomes \({\mathbb{R}}^{n}\) , we have

\begin{theorem}
	Let \(u\) be a positive solution of \eqref{vle} in \({\mathbb{R}}^{n}\) with \(\alpha  \in  \left( {-\infty ,{p}_{S}\left( n\right) }\right)\). Then we have
	\begin{equation}\label{512}
			\mathop{\sup }\limits_{{\mathbb{R}}^{n}}\left( {\frac{{\left| \nabla u\right| }^{2}}{{u}^{2}} + {u}^{\alpha  - 1}}\right)  \leq  C\left( {n,\alpha }\right) \mathop{\sup }\limits_{{\mathbb{R}}^{n}}\left( {{\left| V\right| }^{2} + \left| {\nabla V}\right| }\right).
	\end{equation}
\end{theorem}

\subsection{An example and sharpness of estimate}
The following example illustrates the optimality of estimate \eqref{512}.
	\begin{example}[Entire solution of \eqref{vle}] For any \(\alpha  \in  \left( {1,\frac{n + 2}{n - 2}}\right) \left( {n \geq  3}\right)\) , we take vector field \(V\) as
	
	\[
	V = V\left( {\mu ,n,\alpha }\right) \left( x\right)  = \frac{\Gamma \left( {n,\alpha }\right) }{{\mu }^{2} + {\left| x\right| }^{2}}x,\qquad\forall x \in  {\mathbb{R}}^{n},
	\]
	where \(\mu  > 0\) and
	
	\[
	\Gamma \left( {n,\alpha }\right)  = \frac{4}{\alpha  - 1} - n + 2. 
	\]
	Then the equation \eqref{vle} admits the following global positive solutions:

	\[
	u\left( x\right)  = {\left( \frac{\mu }{{\mu }^{2} + {\left| x\right| }^{2}}\sqrt{\frac{4n}{\alpha  - 1}}\right) }^{\frac{2}{\alpha  - 1}}. 
	\]
	So \({\begin{Vmatrix}{u}^{\alpha  - 1}\end{Vmatrix}}_{{L}^{\infty }\left( {\mathbb{R}}^{n}\right) } = \frac{4n}{\alpha  - 1}{\mu }^{-2}\). Direct computation shows that
	
	\[
	\mathop{\sup }\limits_{{\mathbb{R}}^{n}}\left( {{\left| V\right| }^{2} + \left| {\nabla V}\right| }\right)  = C\left( {n,\alpha }\right) {\mu }^{-2},
	\]
	hence the estimate \eqref{512} is sharp.
\end{example}

\end{document}